\documentclass[11pt,fleqn]{scrartcl}

\usepackage[english]{babel}
\usepackage{amsmath,amsfonts, amsthm,xcolor,amssymb,comment}
\usepackage{paralist}
\numberwithin{equation}{section}
\usepackage{mathtools, mathabx}

\usepackage{hyperref}
\hypersetup{linktoc=none, bookmarksnumbered, colorlinks=true, linkcolor=blue, citecolor=blue}
\usepackage[nameinlink,noabbrev]{cleveref}
\crefformat{equation}{(#2#1#3)}
\crefrangeformat{equation}{(#3#1#4) to~(#5#2#6)}
\crefmultiformat{equation}{(#2#1#3)}{ and~(#2#1#3)}{, (#2#1#3)}{ and~(#2#1#3)}
\crefname{subsection}{subsection}{subsections}
\Crefname{subsection}{Subsection}{Subsections}

\usepackage{xcolor}

\usepackage[hyperpageref]{backref}

\newtheorem{thm}{Theorem}[section]

\newtheorem{cor}[thm]{Corollary}
\newtheorem{prop}[thm]{Proposition}

\theoremstyle{definition}
\newtheorem{rem}[thm]{Remark}
\theoremstyle{remark}

\newcommand{\ds}{\displaystyle}

\newcommand{\R}{\mathbb{R}}

\newcommand{\cH}{\mathcal H}

\newcommand{\de}{\partial}

\DeclareMathOperator{\dive}{div}

\DeclareMathOperator{\Qp}{\mathcal Q_p}

\DeclareMathOperator{\arccosh}{arccosh}

\newcommand\restr[2]{{
  \left.\kern-\nulldelimiterspace 
  #1 
  \vphantom{ |} 
  \right|_{#2} 
  }}
{\left\{\begin{array}{@{}l@{}}}{\end{array}\right.}
\patchcmd{\abstract}{\scshape\abstractname}{\textbf{\abstractname}}{}{}
\makeatletter 
\def\@makefnmark{} 
\makeatother 

\begin{document}

\title{P\'olya-type estimates for the first Robin eigenvalue of elliptic operators}
\author{Francesco Della Pietra$^{*}$ 
\thanks{Dipartimento di Matematica e Applicazioni ``R. Caccioppoli'', Universit\`a degli studi di Napoli Federico II, Via Cintia, Monte S. Angelo - 80126 Napoli, Italia.  \newline 
Email: f.dellapietra@unina.it }
}
\date{}
\maketitle
\begin{abstract}
\noindent{\textsc{Abstract.}} The aim of this paper is to obtain optimal estimates for the first Robin eigenvalue of the anisotropic $p$-Laplace operator, namely:
\begin{equation*}
\lambda_F(\beta,\Omega)=
\min_{\psi\in W^{1,p}(\Omega)\setminus\{0\} } \frac{\displaystyle\int_\Omega F(\nabla \psi)^p dx +\beta \ds\int_{\de\Omega}|\psi|^p F(\nu_{\Omega}) d\cH^{N-1} }{\displaystyle\int_\Omega|\psi|^p dx},
\end{equation*}
where $p\in]1,+\infty[$, $\Omega$ is a bounded, convex domain in $\R^{N}$, $\nu_{\Omega}$ is its Euclidean outward normal, $\beta$ is a real number, and $F$ is a sufficiently smooth norm on $\R^{N}$.
We show an upper bound for $\lambda_{F}(\beta,\Omega)$ in terms of the first eigenvalue of a one-dimensional nonlinear problem, which depends on $\beta$ and on the volume and the anisotropic perimeter of $\Omega$, in the spirit of the classical estimates of P\'olya \cite{po61} for the Euclidean Dirichlet Laplacian.

We will also provide a lower bound for the torsional rigidity
\[
\tau_p(\beta,\Omega)^{p-1} = \max_{\substack{\psi\in
    W^{1,p}(\Omega)\setminus\{0\}}}
\dfrac{\left(\displaystyle\int_\Omega |\psi| \, 
    dx\right)^p}{\displaystyle\int_\Omega F(\nabla\psi)^p dx+\beta \ds\int_{\de\Omega}|\psi|^p F(\nu_{\Omega}) d\cH^{N-1} },
\]
when $\beta>0$. The obtained results are new also in the case of the classical Euclidean Laplacian.
\\[.5cm]
\noindent \textbf{MSC 2020:} 35J25 - 35P15 - 47J10 - 47J30. \\[.2cm]
\textbf{Key words and phrases:}  Nonlinear eigenvalue problems; Robin boundary conditions; Finsler norm; Optimal estimates.
\end{abstract}

\section{Introduction}
The purpose of this work is to provide optimal upper bounds for the first Robin eigenvalue of a Finsler $p$-Laplacian operator with Robin boundary conditions, as well as lower bounds for the corresponding torsional rigidity, in convex domains. This will be achieved by formulating these quantities in terms of a one-dimensional eigenvalue problem. The results presented here are novel, even in the well-established context of the Euclidean Robin-Laplacian operator.
To be more specific, let us first define
\begin{equation}
\label{eigint}
\lambda_{F}(\beta,\Omega)=
\min_{\psi\in W^{1,p}(\Omega)\setminus\{0\} } \frac{\displaystyle\int_\Omega F(\nabla \psi)^p dx +\beta \ds\int_{\de\Omega}|\psi|^pF(\nu_{\Omega}) d\cH^{N-1} }{\displaystyle\int_\Omega|\psi|^p dx},
\end{equation}
where $p\in]1,+\infty[$, $\Omega$ is a bounded Lipschitz domain in $\R^{N}$, $\nu_{\Omega}$ is its Euclidean outward normal, $\beta$ is a real number, $F : \mathbb{R}^N \to [0, +\infty[$, $N\ge 2$, is a convex, even,  $1$-homogeneous and $C^{2}(\mathbb{R}^N\setminus \{0\})$ function such that $[F^{2}]_{\xi\xi}$ is positive definite in $\R^{N}\setminus\{0\}$ (we refer the reader to  \Cref{notation} for all the definitions, here and below used, related the Finsler metric $F$). When $\beta=0$, $\lambda_F(\beta,\Omega)=0$ corresponds to the first (trivial) Neumann eigenvalue, while, when $\beta$ goes to $+\infty$, it reduces to the first Dirichlet eigenvalue $\lambda_{F}^{D}(\Omega)$. Let us observe that, by using $\psi=1$ and $\psi=u_{D}$, the first Dirichlet eigenfunction, in \eqref{eigint}, it holds that
\[
\lambda_{F}(\beta,\Omega)\le \min\left\{ \beta\frac{P_{F}(\Omega)}{|\Omega|},\lambda_{F}^{D}(\Omega)\right\},
\]
where $P_{F}(\Omega)$ is the anisotropic perimeter of $\Omega$ and $|\Omega|$ is its Lebesgue measure. We get in particular that $\lambda_{F}(\beta,\Omega)\to -\infty$ as $\beta\to-\infty$.

If $u\in W^{1,p}(\Omega)$ is a minimizer of \cref{eigint}, then it satisfies
\begin{equation}
    \label{pb_Robin1}
\begin{cases}
-\mathcal Q_p u =\lambda_{F}(\beta,\Omega)\left|u\right|^{p-2}u \quad &\text{in}\ \Omega,\\
F(\nabla u)^{p-1}F_{\xi}(\nabla u)\cdot \nu_{\Omega}+\beta F(\nu_{\Omega}) |u|^{p-2}u= 0 & \text{on}\ \partial\Omega,
\end{cases}
\end{equation}
where 
\[
\mathcal Q_p u=\dive(F^{p-1}(\nabla u)F_\xi(\nabla u)).
\]
In the simplest case where $F(\xi)=\mathcal E(\xi):=\sqrt{\sum_{i}|\xi_{i}|^{2}}$ is the Euclidean norm, then $\mathcal Q_{p}=\Delta_{p}u=\dive(|\nabla u|^{p-2}\nabla u)$ is the 
standard $p$-Laplacian.

In the Dirichlet-Laplacian case, a sharp estimate of the first eigenvalue in terms of perimeter and volume can be given. Indeed
\begin{equation}
\label{polyaclassica}
 \lambda_F^{D}(\Omega) \le (p-1) \left( \frac{\pi_p}{2} \right)^p
\frac{P_F(\Omega)^p}{|\Omega|^p} 
\end{equation}
where $\Omega$ is a bounded convex open set of $\mathbb{R}^N$ and $|\Omega|$ and $P_{F}(\Omega)$ denote respectively the Lebesgue measure and the anisotropic perimeter of $\Omega$, and
\begin{equation}
\label{defpp}
\pi_{p}=2\int_{0}^{1}(1-t^{p})^{-\frac{1}{p}}dt=\frac{2\pi}{p\sin\frac\pi p}.
\end{equation}
The estimate is optimal, being achieved asymptotically when the convex domain $\Omega$ goes to an infinite $N-$dimensional slab $]-a,a[\times \R^{N-1}$ (up to a suitable rotation). 

In the planar, Euclidean, Laplacian case ($N=p=2$ and $F=\mathcal E$), this is a classical result obtained  by P\'olya in \cite{po61}. Then it was generalized, in any dimension, for any $p\in]1,+\infty[$ and for any sufficiently smooth norm $F$ in \cite{mana}. Being $\lambda_{F}(\beta,\Omega)\le \lambda_F^{D}(\Omega)$, the estimate \cref{polyaclassica} is true also in the case of the first Robin eigenvalue, but it is no longer optimal. 
The main contribution of this paper lies in demonstrating that it is possible to refine the P\'olya estimate when switching from Dirichlet to Robin boundary conditions.

\begin{thm}\label{polythm}
    Let $\Omega$ be a bounded, convex, open set of $\R^N$. Then,
    \begin{equation}
      \label{thesis}
\lambda_F(\beta,\Omega) \le \mu_{1}(\beta,s_{0})
\end{equation}
where 
\[
\mu_{1}(\beta,s_{0})=\min_{h\in \mathcal A} \frac{\ds \int_0^{s_{0}} [-h'(s)]^p
  ds+ \beta  h\left(s_{0}\right)^{p}}{\ds\int_0^{s_{0}} [h(s)]^p ds}
\]
where $s_{0}=\frac{|\Omega|}{P_{F}(\Omega)}$, and $\mathcal A$ is the class of positive decreasing functions in $W^{1,p}(0,s_{0})$. The estimate is optimal, being achieved asymptotically when $\Omega$ goes to a suitable slab. 
\end{thm}
\begin{rem}
We draw attention to the fact that this result is new even for the Euclidean norm and for the case of $p=2$.
\end{rem}
\begin{rem}
The minimum $\mu_{1}$ is the first eigenvalue of the nonlinear one dimensional problem
\begin{equation}
\label{1dimx}
    \begin{cases}
    (|X'|^{p-2}X')'+\mu |X|^{p-2}X=0 \quad\text{in}\ (0,s_0),\\
    X'(0)=0, \\
    |X'(s_{0})|^{p-2}X'(s_{0})+ \beta |X(s_{0})|^{p-2}X(s_{0})=0.
    \end{cases}
\end{equation}
This value can be written in terms of the generalized trigonometric functions (when $\beta>0$) or generalized hyperbolic functions (when $\beta<0$) (see \Cref{1D_subsec}). In particular:
\begin{itemize}
\item[\textcolor{gray}{$\bullet$}] if $\beta>0$ 
it holds that (see \Cref{mu_thm} below)
\[
\mu_{1}(\beta,s_{0}) \le 
(p-1) \left( \frac{\pi_p}{2} \right)^p \frac{P_F(\Omega)^p}{|\Omega|^p}.
\]
\item[\textcolor{gray}{$\bullet$}] If $\beta<0$, it is well-known that (see, for example, \cite{gs1,KP,dppstimebasso}) 
\begin{equation}
\label{oldbetaneg}
\lambda_{1}(\beta,\Omega)\le -(p-1)|\beta|^{p'}.
\end{equation}
Also in this case the estimate in \cref{thesis} is sharper. Indeed, the analogous one-dimensional inequality  $\mu_{1}(\beta,s_{0})\le -(p-1)|\beta|^{p'}$ (see \Cref{mu_thm} below) and \cref{thesis} give \cref{oldbetaneg}.
\end{itemize}
\end{rem}
\begin{rem}
In the case $p=2$ (hence in the Laplacian case if $F$ is the Euclidean norm), we get a fine estimate for $\mu_{1}$ by using the expression of $\mu_{1}$ and a Becker-Stark inequality \cite{beckerstark}. So we have the following
\begin{cor}
\label{cor1}
 Suppose that $p=2$, and that the hypotheses of \Cref{polythm} are satisfied, with $\beta>0$. Then
\[
\lambda_{F}(\beta,\Omega) \le \frac{\pi^{2}}{4} \frac{P_{F}^{2}(\Omega)}{|\Omega|^{2}} \frac{1}{1+\frac{2P_{F}(\Omega)}{\beta|\Omega|}}.
\]
\end{cor}
\end{rem}


We observe that other kind of optimal lower bounds, obtained in terms of different geometrical quantities related to $\Omega$ (as the inradius), can be found for example in \cite{Sp92,S,LW} (in the Euclidean case) and in \cite{dppstimebasso}.

The second aim of the paper is also to provide a lower bound to the anisotropic $p$-torsional rigidity with Robin boundary conditions ($\beta>0$), namely the value $\tau_{F}(\beta,\Omega)$ such that  
\begin{equation}
  \label{tors}
\tau_{F}(\beta,\Omega)^{p-1} = \max_{\substack{\psi\in
    W^{1,p}(\Omega)\setminus\{0\} }}
\dfrac{\left(\displaystyle\int_\Omega |\psi| \, 
    dx\right)^p}{\displaystyle\int_\Omega F(\nabla \psi)^p dx+\beta\int_{\de\Omega}|\psi|^{p}F(\nu)d\mathcal H^{N-1}},
\end{equation}
or, equivalently
\begin{equation*}
  \tau_{F}(\beta,\Omega)=\int_\Omega F(\nabla u_p)^p dx +\beta \int_{\de\Omega}|u_{p}|^{p}F(\nu)d\mathcal H^{N-1} = \int_\Omega u_p dx,
\end{equation*}
where $u_p\in W^{1,p}(\Omega)$ is the unique solution of
\begin{equation*}
  \left\{
  \begin{array}{ll}
    -\mathcal Q_p u = 1 &\text{in }\Omega,\\
    F(\nabla u)^{p-1}F_{\xi}(\nabla u)\cdot \nu+\beta|u|^{p-2}u=0 &\text{on }\de\Omega.
  \end{array}
  \right.
\end{equation*}
Let us remark that, by choosing $\psi=1$, or $\psi=u_{\Omega}$ the Dirichlet torsion function, it holds that
\[
\tau_{F}(\beta,\Omega)^{p-1} \ge \max\left\{ \frac{|\Omega|^{p}}{\beta P_{F}(\Omega)},\tau_{F}^{D}(\Omega)^{p-1}\right\},
\]
where $\tau_{F}^{D}(\Omega)$ is the torsional rigidity
with Dirichlet boundary conditions.
The obtained lower bound for $\tau_{p}$ is the following.
\begin{thm}\label{teotor}
  Let $1<p<+\infty$, $\beta>0$ and $\Omega$ be a bounded convex open set in $\R^{N}$. Then,
\begin{equation}
\label{stimator}
\tau_{F}(\beta,\Omega) \ge\left( \frac{p-1}{2p-1}|\Omega| +\frac{1}{\beta^{\frac{1}{p-1}}} \right) \left(\frac{|\Omega|}{P_{F}(\Omega)}\right)^{\frac{p}{p-1}}.
  \end{equation}
\end{thm}
When $\beta=+\infty$, we recover the P\'olya estimate contained in \cite{po61} for the Euclidean Laplacian and then in \cite{mana} for the general case.

The structure of the paper is the following. In \Cref{notation} we recall some useful properties of the Finsler norm, as well as some basic definitions of the anisotropic perimeter and of convex analysis. Moreover, in \Cref{1D_subsec} we recall the main properties of the one dimensional nonlinear eigenvalue problem \Cref{1dimx}. Finally, in \Cref{proofs} we give the proof of the main results.

\section{Notation and preliminaries}
\label{notation}

\subsection{The Finsler norm}
Let
\[
\xi\in \R^{N}\mapsto F(\xi)\in [0,+\infty[
\] 
be a convex, even, $1-$homogeneous function, that is a convex function such that
\begin{equation}
\label{eq:omo}
F(t\xi)=|t|F(\xi), \quad t\in \R,\,\xi \in \R^{N}, 
\end{equation}
 and such that
\begin{equation}
\label{eq:lin}
a|\xi| \le F(\xi)\le b|\xi|,\quad \xi \in \R^{N},
\end{equation}
for some constants $0<a\le b$. Moreover, we suppose that $F\in C^{2}(\mathbb{R}^N\setminus \{0\})$, and that
\begin{equation}
\label{strong}
\nabla^{2}_{\xi}[F^{2}](\xi)\text{ is positive definite in }\R^{N}\setminus\{0\}.
\end{equation}

The hypothesis \cref{strong} on $F$ ensures that the operator 
\[
\Qp[u]:= \dive \left(\frac{1}{p}\nabla_{\xi}[F^{p}](\nabla u)\right)
\] 
is elliptic, namely there exists a positive constant $\gamma$ such that
\begin{equation*}
\sum_{i,j=1}^{n}{\nabla^{2}_{\xi_{i}\xi_{j}}[F^{p}](\eta)
  \xi_i\xi_j}\ge
\gamma |\eta|^{p-2} |\xi|^2, 
\end{equation*}
for some positive constant $\gamma$, for any $\eta \in
\R^N\setminus\{0\}$ and for any $\xi\in \R^N$. 

The polar function of $F$ is
\begin{equation}
\label{polardef}
F^o(v)=\sup_{\xi \ne 0} \frac{\langle \xi, v\rangle}{F(\xi)},	\quad v\in \R^{N}. 
\end{equation}
 It holds that $F^o$ is a convex function
which satisfies properties \cref{eq:omo} and
\cref{eq:lin} (with different constants). Furthermore, 
\begin{equation*}
F(v)=\sup_{\xi \ne 0} \frac{\langle \xi, v\rangle}{F^o(\xi)}.
\end{equation*}
From \cref{polardef} it holds that
\begin{equation}\label{prodscal}
\langle \xi, \eta \rangle \le F(\xi) F^{o}(\eta) \qquad \forall \xi, \eta \in \R^{N}.
\end{equation}
The set
\[
\mathcal W = \{  \xi \in \R^N \colon F^o(\xi)< 1 \}
\]
is Wulff shape centered at the origin. We put
$\kappa_N=|\mathcal W|$, where $|\mathcal W|$ denotes the Lebesgue measure
of $\mathcal W$. More generally, we denote with $\mathcal W_r(x_0)$
the set $r\mathcal W+x_0$, that is the Wulff shape centered at $x_0$
with measure $\kappa_Nr^N$, and $\mathcal W_r(0)=\mathcal W_r$.

The anisotropic distance of $x\in\overline\Omega$ to the boundary of a bounded domain $\Omega$ is the function 
\begin{equation*}
d_{F}(x)= \inf_{y\in \de \Omega} F^o(x-y), \quad x\in \overline\Omega.
\end{equation*}

We stress that when $F=|\cdot|$ then $d_F=d_{\mathcal{E}}$, the Euclidean distance function from the boundary.

It is not difficult to prove that $d_{F}$ is a uniform Lipschitz function in $\overline \Omega$ and
\begin{equation*}
  F(\nabla d_F(x))=1 \quad\text{a.e. in }\Omega.
\end{equation*}
Obviously, $d_F\in W_{0}^{1,\infty}(\Omega)$. 
Finally, the anisotropic inradius of $\Omega$ is the quantity
\begin{equation*}
R_{F}(\Omega)=\max \{d_{F}(x),\; x\in\overline\Omega\},
\end{equation*}
that is the radius of the largest Wulff shape $\mathcal W_{r}(x)$ contained in $\Omega$.

Let us finally recall the definition of anisotropic perimeter of
a set $K\subset \R^N$ in $\Omega$:
\[
P_F(K,\Omega) = \sup\left\{ \int_K
  \dive \sigma dx\colon \sigma \in C_0^1(\Omega;\R^N),\;
  F^o(\sigma)\le 1 \right\}.
\]
The following co-area formula for the anisotropic perimeter
\begin{equation}\label{fr}
  \int_\Omega F(\nabla u) dx = \int_{-\infty}^{+\infty} P_F (\{u>s\},\Omega)\, ds,\quad
  \forall u\in W^{1,1}(\Omega)
\end{equation}
holds, moreover
\[
P_F(K;\Omega)= \int_{\Omega\cap \partial^*K} F(\nu_K) d\mathcal H^{N-1}
\]
where $\mathcal H^{N-1}$ is the $(N-1)-$dimensional Hausdorff measure in $\mathbb R^N$, $\partial^*K$ is the reduced boundary of $F$ and $\nu_F$ is the outer normal to $F$. As usual, we will denote by $P_{F}(K)$ the perimeter in $\R^{N}$, that is $P_{F}(K,\R^{N})$.

\section{A one dimensional $p$-Laplacian eigenvalue problem}\label{1D_subsec}
Here we briefly summarize the definitions and some properties of the $p$-trigonometric functions. These functions are generalizations of the standard trigonometric functions, and they coincide with the standard trigonometric functions when $p=2$. We refer the reader, for example, to \cite{LE,L}.


The function $\arccos_p:[0, 1]\to\R$ is defined as
\begin{equation*}
\arccos_{p}(x)=\ds\int_x^1\frac{dt}{\ds\left(1-{t^p}\right)^\frac 1p}.
\end{equation*}
If $z(t)$ is the inverse function of $\arccos_{p}$, which is defined on the interval $\left[0,\frac{\pi_p}{2}\right]$, where
\[
\pi_p=2\ds\int_{ 0}^{1}\frac{dt}{\ds\left(1-{t^p}\right)^\frac 1p}
= \frac{2\pi}{p\sin\frac{\pi}{p}},
\]
then, the $p$-cosine function $\cos_p$ is the even function defined as the periodic extension of $z(t)$:
\[
\cos_p(t)=\left\{ \begin{split}
& z(t) &&\text{if}\ \ t\in\left[0,\frac{\pi_p}{2}\right],\\
& -z(\pi_p-t) &&\text{if}\ \ t\in\left[\frac{\pi_p}{2}, \pi_p\right], \\
& \cos_p(-t)\  &&\text{if}\ \ t\in\left[-\pi_p, 0\right], \\
\end{split}
\right.
\]
and extended periodically to all $\R$, with period $2\pi_p$; the extension is continuosly differentiable on $\R$. If $p=2$, then $\cos_{p}x$ and $\arccos_{p}x$ coincides with the standard trigonometric functions $\cos x$ and $\arccos x$.

Let now also recall the definitions of the generalized hyperbolic cosine and arccosine functions. The function $\arccosh_{p}$ is defined as
\begin{equation*}
\arccosh_p(x)=\ds\int_1^x\frac{1}{(t^p-1)^\frac1p}dt,\ x\in [1,+\infty[.
\end{equation*}
Its inverse function will be denoted by $\cosh_{p} \colon t\in [0,+\infty[\mapsto[1,+\infty[$. This functio is strictly increasing in $[0,+\infty[$; it can be extended on all $\R$ as $\cosh_{p}(-t)=\cosh_{p}(t)$, $t>0$.
If $p=2$, $\cosh_{p}$ and $\arccosh_{p}$ are the standard hyperbolic functions.

Now we consider the following eigenvalue problem in the unknown $X=X(s)$:
\begin{equation}
\label{1dim}
    \begin{cases}
    (|X'|^{p-2}X')'+\mu |X|^{p-2}X=0 \quad\text{in}\ (0,s_0),\\
    X'(0)=0, \\
    |X'(s_{0})|^{p-2}X'(s_{0})+\beta |X(s_{0})|^{p-2}X(s_{0})=0, 
    \end{cases}
\end{equation}
where $s_{0}$ is a given positive number.

The following result holds (see for example \cite{dppstimebasso}).

\begin{thm}\label{mu_thm}
Let $1<p<+\infty$, $s_{0}>0$ and $\beta \in \R$. Then there exists the smallest eigenvalue $\mu$ of  \cref{1dim}, which has the following variational characterization:
\[
\mu_{1}(\beta,s_0)=\inf_{v\in W^{1,p}(0,s_{0})} \frac{\int_{0}^{s_{0}}\left|v'(s)\right|^{p}ds+\beta v(s_{0})^{p}}{\int_{0}^{s_{0}}\left|v(s)\right|^{p}ds}.
\]
 Moreover, the corresponding eigenfunctions are unique up to a multiplicative constant and have constant sign. The first eigenvalue $\mu_{1}(\beta,s_0)$ has the sign of $\beta$. 
 
 In the case $\beta>0$, the first eigenfunction is
\begin{equation*}
X(s)=\cos_p\left(\left(\frac{\mu_{1}(\beta,s_0)}{p-1}\right)^\frac{1}{p}s \right),\quad s\in (0,s_{0});
\end{equation*}
the eigenvalue $\mu_{1}(\beta,s_0)$ is the first positive value that satisfies
\begin{equation*}
\frac{\mu}{p-1}=\frac{\beta^{p'}}{\cos_p^{-p}\left(\left(\frac{\mu}{p-1}\right)^\frac 1p s_{0} \right)-1},
\end{equation*}
and it holds that
\begin{equation}
\label{stimamubetapos}
\mu_{1}(\beta,s_0) < (p-1)\left(\frac{\pi_p}{2s_{0}}\right)^{p}.
\end{equation}
If $\beta<0$, the first eigenfunction is
\begin{equation*}
X(s)=\cosh_p\left(\left(\frac{-\mu_{1}(\beta,s_0)}{p-1}\right)^\frac{1}{p}s \right),\quad s\in (0,s_{0})
\end{equation*}
the eigenvalue $\mu_{1}(\beta,s_0)$ is the unique negative value that satisfies
\begin{equation*}
-\frac{\mu}{p-1} =\frac{|\beta|^{p'}}{1-\cosh_p^{-p}\left(\left(\frac{-\mu}{p-1}\right)^\frac 1p s_{0} \right)},
\end{equation*}
and it holds that
\begin{equation}
\label{stimamubetaneg}
\mu_{1}(\beta,s_0) \le -(p-1)|\beta|^{p'}.
\end{equation}
\end{thm}

\section{Proof of the main results}
\label{proofs}

\begin{proof}[Proof of \Cref{polythm}]
  Let $g(t)=g(d_F(x))$, where $g$ is a nonnegative, increasing,
  sufficiently smooth function, where $d_{F}(x)$ is the distance of $x\in \Omega$ from $\de\Omega$. 
  Being  $F(\nabla d_F)=1$, by the coarea formula it holds that
 \[
  \int_\Omega F(\nabla g(d_F(x)))^p dx = \int_0^{R_F(\Omega)} g'(t)^p dt \int_{\{d_F=t\}}
  \frac{1}{|\nabla d_F|} d\mathcal H^{N-1} = \int_0^{R_F(\Omega)} g'(t)^p
  P(t)dt,
  \]
  where $P(t)=P_F(\{x\in \Omega\colon d_F>t\})$, and $R_F(\Omega)=\sup_{x\in\Omega} d_{F}(x)$; similarly,
  \[
  \int_\Omega g(d_F(x))^p dx = \int_0^{R_F(\Omega)} g(t)^p P(t)dt, \quad \int_{\de\Omega} g(d_{F}(x))F(\nu) d\mathcal H^{N-1}= g(0) P_{F}(\Omega).
  \]
   Using $g(d_F(x))$ as test function in the Rayleigh quotient
  of \cref{pb_Robin1}, we have
  \begin{equation}
    \label{polyacarvar}
  \lambda_F(\beta,\Omega) \le \frac{\ds \int_0^{R_F(\Omega)} g'(t)^p
    P(t)\,dt + \beta P_{F}(\Omega)g(0)^{p}}{\ds \int_0^{R_F(\Omega)} g(t)^p P(t)\,dt}.
  \end{equation}
  Now, we perform the change of variable
  \[
  s =  \frac{C}{|\Omega|} A(t),
  \]
  where we have denoted by $A(t)=|\{x\in
  \Omega\colon d_F(x)>t\}|$, and $C$ is a positive constant which will
  be chosen in the next. Observe that $A(0)=|\Omega|$. Let $h(s)$ be the function such that
  \[
  h(s)=g(t).
  \]
  We stress that $h(s(t))$ is decreasing. Then, substituting $h$ in
  \cref{polyacarvar}, it follows that
  \begin{equation*}
  \lambda_F(\beta,\Omega) \le \frac{C^p}{|\Omega|^p}
  \frac{\displaystyle\int_0^{C} \left[-h'\left(
      \frac{C}{|\Omega|}A(t)\right)\right]^p P(t)^p
    \left[\frac{C}{|\Omega|} P(t) \right] dt
    + 
    \beta \frac{|\Omega|^{p-1}}{C^{p-1}} P_{F}(\Omega) h(C)^{p}
    }
  {\displaystyle\int_0^{R_F(\Omega)} \left[ h \left(
        \frac{C}{|\Omega|}A(t)\right)\right]^p
  \left[\frac{C}{|\Omega|} P(t) \right] dt}
\end{equation*}
and then, by the monotonicity of $P(t)$ we get
\begin{equation}
  \label{fine}
  \lambda_F(\beta,\Omega) \le
  \frac{P_F(\Omega)^p}{|\Omega|^p} C^p
  \frac{\ds	\int_0^{C} [-h'(s)]^p ds+\beta \frac{|\Omega|^{p-1}}{C^{p-1} P_{F}(\Omega)^{p-1}} h(C)^{p}}{\ds\int_0^{C} [h(s)]^p ds}.
\end{equation}
We want to minimize the quantity in the right-hand side of equation \cref{fine}. As matter of fact, such minimum
does not depend on $C$ so we can set $C=s_{0}=\frac{|\Omega|}{P_{F}(\Omega)}$. This gives us the following expression:
\begin{equation}
  \label{minpb}
  \lambda_F(\beta,\Omega) \le  \mu_{1}(\beta,s_{0}), 
\end{equation}
where 
\[
\mu_{1}(\beta,s_{0})=\min_{\mathcal A} \frac{\ds \int_0^{s_{0}} [-h'(s)]^p
  ds+\beta h\left(s_{0}\right)^{p}}{\ds\int_0^{s_{0}} [h(s)]^p ds}
\]
and $\mathcal A$ is the class of positive decreasing functions $h\in
W^{1,p}(0,s_{0})$. 

The optimality of the inequality follows by using a similar argument of \cite[Proposition 5.1]{dppstimebasso}. More precisely, what can be proved is the following.

\begin{prop} Let $\Omega_\ell=]-\frac a2,\frac a2[\times]-\frac \ell 2,\frac \ell 2[^{N-1}$. Then, up to a suitable rotation of $\Omega_{\ell}$, it holds that
\[
\lim_{\ell\to +\infty}\frac{\lambda_1(\beta,\Omega_\ell)}{\mu_\ell}=1
\]
where $\mu_{\ell}$ is the first eigenvalue of \eqref{1dim} in $(0,s_{\ell})$, with $s_{\ell}=\frac{|\Omega_{\ell}|}{P_{F}(\Omega_{\ell})}$.
\end{prop}
\end{proof}


\begin{proof}[Proof of \Cref{cor1}]
The one-dimensional eigenvalue in this case satisfies
\[
\tan\left(\sqrt\mu_{1}s_{0}\right) = \frac{\beta}{\sqrt{\mu_{1}}},
\]
where $s_{0}=\frac{|\Omega|}{P_{F}(\Omega)}$.
Recalling the Becker-Stark \cite{beckerstark} inequality
\begin{equation}
\label{beckerstark}
\frac{2t}{\frac{\pi^{2}}{4}-t^{2}} \le \tan t, 
\qquad 0\le t < \frac{\pi}{2};
\end{equation}
 it holds that
\[
\frac{\beta}{\sqrt{\mu_{1}}} \ge \frac{2s_{0}\sqrt \mu_{1}}{\frac{\pi^{2}}{4}-s_{0}^{2}\mu_{1}}.
\]
Hence, rearranging and recalling \Cref{polythm} we get
\[
\lambda_{1}(\beta,\Omega) \le 
\mu_{1} \le \frac{\pi^{2}}{4} \frac{P_{F}^{2}(\Omega)}{|\Omega|^{2}} \frac{1}{1+\frac{2P_{F}(\Omega)}{\beta|\Omega|}}.
\]
\end{proof}


\begin{proof}[Proof of \Cref{teotor}]
 We use the same notation and argument of
 the proof of \Cref{polythm}. For a test function
 $g(t)=g(d_F(x))$, with $g$ nonnegative increasing sufficiently smooth
 function, we have that
 \begin{multline*}
   \label{eq:5}
  \tau_{F}(\beta,\Omega)^{p-1} \ge 
  \dfrac{\left(\displaystyle\int_0^{R_F(\Omega)}g(t)
      P_F(t)dt\right)^p}{\displaystyle\int_0^{R_F(\Omega)} g'(t)^p
    P_F(t)dt +\beta g(0)^{p }P_{F}(\Omega)^{p}}
    = \\ = \dfrac{\left(\displaystyle\int_0^{R_F(\Omega)}g'(t)
      A(t)dt+g(0)|\Omega|\right)^p}{\displaystyle\int_0^{R_F(\Omega)} g'(t)^p
    P_F(t)dt+\beta g(0)^{p }P_{F}(\Omega)^{p}}.  
 \end{multline*}
  Last equality follows performing an integration by parts. Substituting
  \[
  g(t)=\int_0^{t} \left(\frac{A(s)}{P_F(s)}\right)^{1/(p-1)}ds + c,
  \]
  where $c$ is a constant to be chosen later, it follows that
  \begin{equation}
  \tau_{F}(\beta,\Omega)^{p-1} \ge 
  \dfrac{\left(\ds\int_0^{R_F(\Omega)}
  \frac{A(t)^{\frac{p}{p-1}}}{P(t)^{\frac{1}{p-1}}} dt+c|\Omega|\right)^{p}}{\ds\int_0^{R_F(\Omega)}
  \frac{A(t)^{\frac{p}{p-1}}}{P(t)^{\frac{1}{p-1}}} dt+\beta c^{p }P_{F}(\Omega)^{p}}.
  \end{equation}
  Maximizing in $c$, it holds that the optimal choice is $c=\frac{|\Omega|^{\frac{1}{p-1}}}{\beta^{\frac{1}{p-1}}P_{F}(\Omega)^{\frac{p}{p-1}}}$ and then
  \begin{equation}
  \label{catena2}
  \tau_{F}(\beta,\Omega) \ge \ds\int_0^{R_F(\Omega)}
  \frac{A(t)^{\frac{p}{p-1}}}{P(t)^{\frac{1}{p-1}}} dt + \frac{1}{\beta^{\frac{1}{p-1}}} \left(\frac{|\Omega|}{P_{F}(\Omega)}\right)^{\frac{p}{p-1}}.
  \end{equation}
  Now, let us observe that 
  \begin{multline}\label{catena}
  \ds\int_0^{R_F(\Omega)}
  \frac{A(t)^{\frac{p}{p-1}}}{P(t)^{\frac{1}{p-1}}} dt
=  \int_0^{R_F(\Omega)}
  \frac{A(t)^{\frac{p}{p-1}}[-A'(t)]}{P(t)^{\frac{p}{p-1}}} dt
 = \\ = \frac{p-1}{2p-1}
  \frac{|\Omega|^{\frac{2p-1}{p-1}}}{P_F(\Omega)^{\frac{p}{p-1}}} +
\frac{p}{2p-1}
\int_0^{R_F(\Omega)}
\left(\frac{A(t)}{P(t)}\right)^{\frac{2p-1}{p-1}} [-P'(t)] dt \ge \\
\ge
\frac{p-1}{2p-1} \frac{|\Omega|^{\frac{2p-1}{p-1}}}{P_F(\Omega)^{\frac{p}{p-1}}}.
\end{multline}

Given that $P(t)$ is decreasing, we explicitly state that the second equality in \cref{catena} is followed by an integration by parts, observing that for
\[
A(t) =\int_t^{R_F(\Omega)} P(t) dt \le P(t)(R_F(\Omega)-t), 
\]
we have
\[
\frac{|A(t)|^{\frac{2p-1}{p-1}}}{P(t)^{\frac{p}{p-1}}}\rightarrow
0\text{ as }t\rightarrow R_F(\Omega).
\]
In conclusion, by \cref{catena} and \cref{catena2} we have
\[
\tau_{F}(\beta,\Omega) \ge\left( \frac{p-1}{2p-1}|\Omega| +\frac{1}{\beta^{\frac{1}{p-1}}} \right) \left(\frac{|\Omega|}{P_{F}(\Omega)}\right)^{\frac{p}{p-1}}
\]
and this conclude the proof.
\end{proof}

\section*{Acknowledgements}
This work has been partially supported by the
PRIN PNRR 2022 ``Linear and Nonlinear PDE’s: New directions and Applications'', by GNAMPA of INdAM, by  the FRA Project (Compagnia di San Paolo and Universit\`a degli studi di Napoli Federico II) \verb|000022--ALTRI_CDA_75_2021_FRA_PASSARELLI|.

{\small

\bibliographystyle
{plain}


\end{document}